\documentclass[final,letterpaper,11pt]{article}
\usepackage{amsfonts}
\usepackage{graphicx}

\include{tcilatex}

\def\vert{\mathop{\textstyle vert}}
\def\Conv{\mathop{\textstyle conv}}

\newcommand\x{\mathbf{x}}
\renewcommand\l{\mathbf{l}}

\begin{document}

\author{Robert Erdahl, Konstantin Rybnikov}
\title{Supertopes}
\date{8 October 2002}

\maketitle

\abstract{This document reflects the state of the subject as of the end of 2002. It is
an updated and \emph{corrected} version of a preprint initially produced in the end of
the summer of 2001 by Rybnikov at the Department of Mathematics of Cornell University.
}

\section{\protect\large \ Introduction}

A lattice Delaunay cell is \emph{perfect} (a.k.a. extreme) if its Delaunay sphere is
the only ellipsoid that circumscribes this cell. Such Delaunay cells correspond,
roughly in one-to-one fashion, to extreme hypermetrics, initially studied in analysis
and combinatorics (Deza and Laurent 1997). When the lattice is affinely mapped to
$\mathbb{Z}^{n}$, the perfect Delaunay cell is mapped to a $\mathbb{Z}^{n}$-polytope,
which is circumscribed with an empty \emph{perfect ellipsoid }$f(\mathbf{x})=c$: the
defining property
of such ellipsoid is that the \emph{inhomogeneous }quadratic form $f(\mathbf{%
x})$ can be reconstructed in a unique way from its minimum on $\mathbb{Z}%
^{n} $ and all representations of this minimum. A perfect ellipsoid is an exact
inhomogeneous analog of the notion of \emph{perfect form,} introduced by Korkin and
Zolotareff (1873) and later studied by Voronoi (1908-1909), Coxeter (1951), Conway,
Sloane (1988), Martinet (1996) etc. It might not be so well known, but Korkin \&
Zolotareff (1873) came to the discovery of important lattices (and forms)
$E_{n},A_{n},D_{n}$ by trying to construct infinite series of perfect forms! Point
lattices are very important to many areas of algebra, number theory, geometry,
combinatorics, cryptography, communication theory, and the theory of approximations:
e.g., see Conway and Sloane (1999). To understand properties of lattices we often need
to understand certain polytopes
associated with these lattices. Delaunay polytopes (also called \textit{holes%
}) form one of the most important classes of such polytopes.

Positive definite quadratic forms (PDQFs) in $n$ variables make an \emph{%
open }cone $\mathfrak{P}_{n}$ of dimension $N=\frac{n(n+1)}{2}$ in $Sym_{n}(%
\mathbb{R)\cong R}^{N},$ the space of quadratic forms, or symmetric matrices. The
boundary of $\mathfrak{P}_{n}$ consists of positive semi-definite quadratic forms
(referred to as PQFs). PDQFs serve as algebraic representations of \textit{point
lattices}. There is a natural one-to-one correspondence between isometry classes of
$n$-dimensional lattices and
integral equivalence classes (i.e. with respect to substitutions $\mathbf{x}%
^{\prime }\mathbf{=}A\mathbf{x,}$ $A\in GL_{n}(\mathbb{Z})$) of PDQFs in $n$
variables.

Conjugation by a fixed matrix from $GL_{n}(\mathbb{Z})$ is an invertible
linear operator on $Sym_{n}(\mathbb{R)}$. Therefore, conjugation \ defines a
homomorphism $\mathcal{V}$ from $GL_{n}(\mathbb{Z})$ to $GL_{N}(\mathbb{Z}),$
and $GL_{n}(\mathbb{Z})$ acts pointwise on $Sym_{n}(\mathbb{R)}$. Two
subsets of $Sym_{n}(\mathbb{R})$ are called arithmetically equivalent if
they are equivalent with respect to the action of $\mathcal{V}(GL_{n}(%
\mathbb{Z}))$.

\begin{definition}
\emph{A partition} $\mathfrak{R}$ of $\mathfrak{P}_{n}$ \emph{into }relatively open
convex polyhedral \emph{cones with apex at 0 is called a\ reduction partition%
} \emph{if: (1) it is invariant with respect to }$GL_{n}(\mathbb{Z})$\emph{;
(2) there are only finitely many arithmetically non-equivalent cones in this
partition; (3 ) for each cone }$C$\emph{\ of }$R$\emph{\ and any PQF }$%
\varphi $\emph{\ in }$n$\emph{\ indeterminates, }$\varphi $\emph{\ can be }$%
GL_{n}(\mathbb{Z})$\emph{-equivalent to at most finitely many forms from }$C$%
\emph{.}
\end{definition}

Voronoi defined two polyhedral reduction partitions of $\mathfrak{P}_{n}$:\ these are
the tilings by \emph{perfect domains} and \emph{domains for lattice types,} also
called\emph{\ L-domains. }

\begin{definition}
\emph{A convex polyhedron } $P$ \emph{\ in } $\mathbb{R}^{n}$ \emph{\ is called a
Delaunay cell of a lattice } $L$ \emph{\ with respect to a positive quadratic form }
$\varphi (x)$ \emph{\ if: (1) for each face } $F$ \emph{\ of } $P$ \emph{\ we have }
$\Conv(L \cap F)=F$ \emph{; (2) there is a quadric } $Q_{(P,\varphi )}(x)=0,$ \emph{\
circumscribed about } $P$ \emph{, whose quadratic part is } $\varphi (x)$ \emph{; (3)
no points of } $L\backslash P$ \emph{\ satisfy } $Q_{(P,\varphi )}(x)\leq 0.$
\end{definition}

When $\varphi (\mathbf{x)}=\sum x_{i}^{2}$, this definition gives the
classical concept of Delaunay cell in $\mathbb{E}^{n}$ (1924, 1937).
Delaunay tilings can be defined not only for lattices, but for any
reasonable discrete point sets. These tilings have enormous applications in
computational geometry, numerical methods, CAD, the theory of lattices,
mathematical crystallography, etc.

Two $n$-forms $\varphi _{1}$ and $\varphi _{2\text{ }}$ belong to the same $%
L $-type if the Delaunay tilings of $\mathbb{Z}^{n}$ with respect to these
forms are affinely equivalent (the notion of $L$-type is, in fact, due to
Voronoi (1908, vol. 133)). Each $L$-type domain is, of course, the union of
infinitely many convex cones that are equivalent with respect to $GL_{n}(%
\mathbb{Z})$, acting pointwise on $Sym_{n}(\mathbb{R)}$. $L$-type domains form a
\emph{reduction partition} of $\mathfrak{P}_{n}$.

The notions of Delaunay tiling and $L$-type are extremely important in the study of
extremal and group-theoretic properties of lattices. For example, the analysis of
Delaunay cells in the famous Leech lattice conducted by Conway, Sloane and Borcherds
showed that 23 ''deep holes'' (Delaunay cells of radius equal to the covering radius of
the lattice) in the Leech lattice\ are in one-to-one correspondence with even
unimodular 24-dimensional lattices, classified by  Niemeier (1973) that, in turn, give
rise to 23 ''gluing''\ constructions of the Leech lattice from root lattices. Barnes
and Dickson (1967) and Dickson (1968) and, later, Delaunay et al. (1969, 1970) proved
that the closure of each $N$-dimensional $L$-type domain has at most one local minimum
of the ball covering density,
and if such a minimum exists, the group of $GL_{n}(\mathbb{Z})-$%
automorphisms of the domain maps this form to itself. Using this approach
Delaunay et al (1963, 1970) found the best lattice coverings in $\mathbb{E}%
^{4}$ and $\mathbb{E}^{5}.$ The theory of $L$-types has numerous connections to
combinatorics, in particular, to cuts, hypermetrics, and regular graphs (see Deza \&
Laurent 1997), and algebraic geometry. For example, Valery Alexeev (2001, 2002)
recently discovered an interesting connection between $L$-types of $g$-dimensional PQFs
and certain functorial compactification $\overline{AP}_{g}$ of the moduli space of
abelian varieties $A_{g}$.

The $L$-type partition of the cone of PQFs is closely related to the theory
of \emph{perfect forms}. The \emph{arithmetic minimum }of a form $\varphi (%
\mathbf{x})$ is its minimum on $\mathbb{Z}^{n}$. The integral vectors on
which this minimum is attained are called the \emph{minimal vectors} of $%
\varphi $: these vectors have the minimal length among all vectors of $%
\mathbb{Z}^{n}\backslash Ker(\varphi ),$ when $\varphi $ is used as the
metrical form. For each PDF $\varphi $ with the set of minimal vectors $%
\mathsf{P\subset }$ $\mathbb{Z}^{n}$ there is a $\mathsf{P}$-domain $\Pi _{%
\mathsf{P}}$ defined by $\Pi _{\mathsf{P}}=\{\sum_{\mathbf{p\in }\mathsf{P}%
}\omega _{\mathbf{p}}(\mathbf{p\cdot x})^{2}\ |\ \omega _{\mathbf{p}}>0\}.$ The
$\mathsf{P}$-domains form a \emph{reduction partition }of $\mathfrak{P}_{n}$ which is
called the \emph{perfect }partition. Form $\varphi (\mathbf{x})$ is called
\emph{perfect} if it can be reconstructed up to scale from all representations of its
arithmetic minimum. Indeed, uniqueness requires the existence of at least $n(n+1)$
minimal vectors.\ The perfect domains are open polyhedral $N$-dimensional cones, which
fit together facet-to-facet to tile $\mathfrak{P}_{n}$. The $\mathsf{P}$-domains with
lesser dimension are relatively open faces of this tiling. In fact,
$\mathsf{P}$-domains from a partition of the convex hull of all rational rank one
forms:

\begin{equation}
\bigsqcup_{\mathsf{P}}\Pi _{\mathsf{P}} = \Conv \{(\mathbf{p\cdot x}%
)^{2}\ |\ \mathbf{p\in }\mathbb{Q}^{n}\}\subset \overline{\mathfrak{P}}_{n},
\end{equation}

where the union is taken over all sets of minimal vectors for PDQFs in $n$ variables.
Each extreme ray of this tiling lies on $\partial \mathfrak{P}_{n}$.\

Intuitively, perfect lattices are those that have a large supply of minimal
vectors, although a perfect $n$-lattice for $n>8$ is not always spanned by
its minimal vectors. A perfect form $\varphi (\mathbf{x})$ can obviously be
described as a hyperplane in $Sym_{n}(\mathbb{R)}$ that contains $N+1$
integer points whose coordinates are the images of the minimal vectors $\{%
\mathbf{v}_{k}|\,k=1,...,2s\}$ under the Veronese-Voronoi mapping $V:\mathbf{%
v}_{k}\rightarrow \{v_{k}^{i}v_{k}^{j}|\,1\leq i\leq j\leq n\}.$ In fact,
the intersection of the half-spaces (not containing \textbf{0}) defined by
the hyperplanes corresponding to perfect forms is a ''polyhedron'' with
infinitely many faces, called the Voronoi polyhedron. Voronoi showed that
this polyhedron has only finitely many faces which are not $GL_{n}(\mathbb{Z}%
)$-equivalent, and therefore in each dimension there are only finitely many perfect
forms up to $GL_{n}(\mathbb{Z})$-equivalence. Perfect forms play an important role in
lattice sphere packings. Voronoi's theorem (1908, vol. 133) says that a form is\
extreme---i.e., a maximum of the packing density---if and only if it is perfect and
eutactic (see Coxeter, 1951 for details).

A cone in $Sym_{n}(\mathbb{R)}$ spanned by the images, under $V,$ of the minimal
vectors of a perfect form is called a perfect cone. The union of all perfect cones
corresponding to forms integrally equivalent to $\varphi $ is called the perfect domain
of $\varphi $. For each perfect cone there are infinitely many
$GL_{n}(\mathbb{Z})$-equivalent ones, so the perfect domain of $\varphi $ \ consists of
infinitely many equivalent perfect cones, just like an $L$-type domain consists of
infinitely many convex cones.\emph{\ } A fundamental theorem of Voronoi (1908, vol.
133) in the interpretation of Ryshkov (in Delaunay and Ryshkov 1971) says that the cone
of PQFs is tiled face-to-face by perfect cones. Since there are only finitely many
non-equivalent perfect forms, there is a finite set of perfect cones in
$Sym_{n}(\mathbb{R)}$ such that each form in $n$ variables is equivalent to a form from
one of these cones. Voronoi gave an algorithm finding all perfect domains for given
$n$. This algorithm is known as Voronoi's reduction with perfect forms. \ For the
computational analysis of his algorithm and its improvements see Martinet
(1996). Perfect forms have been completely classified in dimensions $n\leq 7$%
, but already for $n=9$ there are billions of them. The theory of perfect forms was
used for finding the best lattice packings in low dimensions and for classifying
maximal finite subgroups of $GL_{n}(\mathbb{Z})$ for small values of $n$ (Ryshkov
1972). Among important recent developments in the theory of perfect forms and $L$-types
is our result (R.E. and K.R. 2002) showing that, to the contrary of the Voronoi (1909)
and Dickson (1972) conjecture, the $L$-partition of the cone of PQFs is \emph{not} a
refinement of the perfect partition for $n\geq 6$. We are now trying to find
connections between $L$ and perfect partitions. This work is based on studying the
pattern of scalar products between the Voronoi vectors of a PDQF $\varphi $ and the
minimal vectors of the perfect form to the domain of which $\varphi $ belongs. One of
distant goals in this direction is to understand the connections between the packing
and covering problems for lattices.

The notion of perfect Delaunay cell is an inhomogeneous analog of the notion
of perfect form.

\begin{definition}
\emph{Let }$P\subset R^{n}$\emph{\ be a Delaunay polytope in }$Z^{n}$\emph{\
for a PDQF }$\varphi (x)$\emph{. }$P$\emph{\ is called perfect if the only
ellipsoid circumscribed about }$P$\emph{\ is the one defined by the form }$%
\varphi (x)$\emph{. In this case the inhomogeneous quadratic form defining
this ellipsoid is referred to as a }perfect ellipsoid.
\end{definition}

The vertices of a perfect Delaunay polytope are analogs of the minimal vectors for a
perfect form: the minimal possible number of vertices of a perfect Delaunay cell is
$\frac{n(n+1)}{2}+n,$ while the minimal number of minimal vectors of a perfect forms is
$\frac{n(n+1)}{2}$. Perfect Delaunay polytopes are important not only to the theory of
lattices, but also to the theory of cuts and hypermetrics started in analysis by I. J.
Schoenberg in 1935-37 (for references see Deza \& Laurent, 1997). In our research we
are primarily interested in geometric aspects of perfect Delaunay cells and ellipsoids,
although we keep in mind possible connections with distance regular graphs: e.g. the
1-skeletons of the Delaunay cells of $E_{6}$ and $E_{7}$ are well-known
strongly-regular graphs---the Schlafli and Gosset graphs.

1-dimensional $L$-types that are interior to the cone of PQFs are very rare
in low dimensions. They first occur in dimension 4: $D_{4}$ has an extreme $%
L $-type. Any new perfect Delaunay cell would give a new example of an extreme
$L$-type. \emph{Prior to our work only finitely many examples of extreme $L$-types and
perfect ellipsoids have been known.} The significance of extreme $L$-types is much due
to their relation to the structure of Delaunay and Dirichlet-Voronoi tilings.
Dirichlet-Voronoi polytopes of higher-dimensional lattices are important to the theory
of quasicrystalls, coding theory, information quantization, etc (e.g. see Conway \&
Sloane (1999), Senechal (1995)). The significance of extreme $L$-types for geometry of
lattices is illustrated by the following theorems.

\begin{theorem}
(Ryshkov 1998; Erdahl 2000) The Voronoi polytope for any PDQF $\varphi $ is the
Minkowski sum of Voronoi polytopes for quadratic forms lying on the extreme rays of the
L-cone of $\varphi $. Their arrangement in space is
determined by $\varphi .$ The Delaunay tiling of $\mathbb{Z}^{n}$ for a PQF $%
\varphi $ is the intersection of the Delaunay tilings of forms lying on the
extreme rays of the $L$-cone of $\varphi $.
\end{theorem}

These propositions are dual formulations of the same fact. For example, the
Delaunay tiling of the plane with respect to the form $%
f=x^{2}+y^{2}+(x-y)^{2}$ is just the intersection of the Delaunay tilings
for forms $x^{2},$ $y^{2}$ and $(x-y)^{2}$ : by vertical strips $k<x<k+1,$ $%
k\in \mathbb{Z}$, \ by horizontal strips $k<y<k+1,$ $k\in \mathbb{Z}$, and
by slanted strips $k<x-y<k+1,$ $k\in \mathbb{Z}.$ The Voronoi cell for $f$
is the sum of three segments orthogonal with respect to $f$ to these three
families of strips. Therefore, Voronoi polytopes for extreme $L$-types are
''building blocks'' for Voronoi polytopes of arbitrary forms, while Delaunay
tilings for extreme rays of the $L$-cone of $f$ are coarsenings of the
Delaunay tiling of $f$.

Very few perfect Delaunay polytopes in low dimensions are known.  It is known that
there are no perfect Delaunay polytopes in dimensions less than 6 (Erdahl 1975). The
canonical examples of perfect Delaunay polytopes are the Gosset polytopes in
$\mathbb{E}^{6}$ and $\mathbb{E}^{7}$, e.g. see Coxeter (1934, 1973), Erdahl (1992),
Deza and Grishukhin (1995, 1996). Other examples include two 16-dimensional polytopes
in $BW_{16}$ and its sublattice, three 15-dimensional polytopes in sections of
$BW_{16}$, two polytopes in 22 and 23 dimensional sections of the Leech lattice (Deza
et al. 1992, 1996; Deza \& Laurent 1997). Most of these examples are manifestations of
such phenomena as extreme sets of equiangular lines and extreme spherical two-distance
sets that have been intensively studied in algebraic combinatorics (see Lemmens and
Seidel, 1973).

\section{Constructing Perfect Delaunay Cells}

In 2001 we constructed an infinite series of big Delaunay cells $%
\Upsilon _{n}$ not coming from equiangular lines or two-distance sets. This series
starts from (the affine image in $\mathbb{Z}^{n}$ of) Gosset polytope in $E_{6}.$ The
construction of this series is based on the\emph{\ infinite series of Delaunay
simplexes of relative volume}$\ n-3$\emph{\ }found by Erdahl and Rybnikov (2002b). This
is the best known infinite series of big Delaunay simplexes; it improves upon Ryshkov's
(1973) series of Delaunay simplexes of volume $r$ in dimension $n=2r+1$ (it is
interesting that Ryshkov's series of big simplexes is also related to his 1973 series
of perfect lattices not generated by its perfect vectors). Polytopes $\Upsilon _{n}$
are constructed by supplementing the vertices of the simplex of volume $n-3$ in some
very special way (see below). This construction generalizes the embedding of the
simplex of volume $n-3$ into the Gosset polytope in $\mathbb{E}^{6}$.  We  hope to
generalize the lattice $E_{8}$ to an infinite series of lattices with interesting
geometric and arithmetic properties.

Lattice Delaunay polytopes of large relative volume and/or many vertices are of special
interest to the study of lattice $L$-types and perfect forms. Polytopes with many
vertices normally occur in highly symmetric lattices, such as $E_{n}$ $(n=6,7),$
Barnes-Wall lattice, $BW_{16}$, sections of Leech lattice, $\Lambda _{24}$, etc.
Simplexes of large volume are found in lattices with special symmetries (e.g. $\Lambda
_{24}$) and their perturbations (e.g. some $L$-type domains with extreme ray of type
$E_{n} \;\; (n=6,7).$ \ Delaunay simplexes are very special cases of empty lattice
simplexes that have been attracting interest of mathematicians due to their importance
in integer programming (e.g. see Haase and Zielger, 2000). While there are pretty sharp
results on empty lattice simplexes, not much is  known about Delaunay empty lattice
simplexes. The volume of an empty lattice simplex can be arbitrary for $n\geq 4$, but
the volume of a Delaunay $n$-simplex is, according to Lov\'{a}sz (unpublished), bounded
from above by $n! (2^{n} / \binom{ 2n}{n})$ (see Deza \& Laurent (1997) for a proof).
It is not even known if the maximal volume of a Delaunay simplex grows linearly,
polynomially, or exponentially in the dimension. The biggest Delaunay simplex, we know
of, lives in the Leech lattice and has volume \ $ 85n$ ($n=24$), but the best infinite
series known is linear: $n-3$. It would be interesting to improve Lov\'{a}sz bound,
since it would improve an upper bound on the number of lattice points one has to check
to verify if a given lattice polytope is Delaunay (Deza \& Laurent, 1997).  Rybnikov is
also interested in \emph{algorithmic approaches to determining whether an empty lattice
polytope is Delaunay.} This problem may be important to cryptography, as it is an
inhomogeneous counterpart of the \emph{shortest vector problem.} He is interested in
both, algorithms for computers, and ''algorithms for humans'', i.e. analytical
techniques, involving arithmetics and geometry, that allow proving the Delaunay
property in the cases when the form is very symmetric, similar to the method of
projective inequalities (Dieter, 1975; Anzin, 1991) for the problem of determining the
shortest vectors.

\subsubsection{Supertopes}

We constructed the following infinite series of \textit{supertopes} $%
\Upsilon _{n}$ in $\mathbb{Z}^{n}$ for all $n\geq 6$. Below are the vertices
of $\Upsilon _{n}$:

\begin{center}
\begin{tabular}{|c|c|c|}
\hline
$\lbrack 0^{n}]\times 1$ & $ [1,0^{n-1};0]\times (n-1) $ & $%
[1,0^{n-2};-1]\times (n-1) $ \\ \hline $[1^{2},0^{n-3};-1]\times \frac{ (n-1)(n-2)}{2}$
&
$ [ 0,1^{n-2};-(n-4)]\times (n-1)$ & $[1^{n-1};-(n-3)]\times 1$ \\
\hline
\end{tabular}
\end{center}

Here we use a short-hand notation for families of vectors obtained from some
$n$-vector by all circular permutations in strings of symbols that are
separated by commas and bordered on the sides by semicolons and/or brackets.
We realized the equation of the ellipsoid circumscribed about $\Upsilon _{n}$
should have the following form:

\QTP{Body Math}
\begin{equation}
d\sum_{\mathbf{i=1}}^{\mathbf{n-1}}x_{\mathbf{i}}^{2}+2m\sum_{\mathbf{1\leq
i<j\leq n}}^{\mathbf{n-1}}x_{\mathbf{i}}x_{\mathbf{j}}+2e\sum_{\mathbf{i=1}%
}^{\mathbf{n-1}}x_{\mathbf{i}}x_{\mathbf{n}}+bx_{\mathbf{n}}^{2}-d\sum_{
\mathbf{i=1}}^{\mathbf{n-1}}x_{\mathbf{i}}+l_{\mathbf{n}}x_{\mathbf{n}}=0
\end{equation}
\pagebreak
\begin{theorem}
\ (Erdahl, Rybnikov, Kemp, Saliola 2001) Equation
$\mathbf{{x}^{t}}\mathbf{Q_n}\mathbf{{x}+} \mathbf{L} \mathbf{{x}=0}$ defines an
ellipsoid circumscribed about polytope $\Upsilon _{n}$ for $n\geq 6$. Here
$\mathbf{Q_n}=(q_{ij})$ is a symmetric positive matrix of the above form where

\begin{center}
$q_{ii}=2-5n+n^{2}$ \emph{for} $i<n$,\end{center}

 \begin{center}
 $q_{nn}=-2-3n+n^{2}$,\end{center}

\begin{center}
 $q_{ij}=12-7n+n^{2}$ \emph{for} $ i<j<n$,\end{center}

\begin{center}
$\ q_{in}=+4-5n+n^{2}$ for $i\neq n$.\end{center}

\noindent $\mathbf{L}$ is a linear functional defined by

\begin{equation}
\l \cdot \x= \sum_{\mathbf{i=1}}^{\mathbf{i=n-1}}
(-q_{\mathbf{ii}})x_{\mathbf{i}}+l_{\mathbf{n}}x_{\mathbf{n}}
 =(-2+5n-n^{2})\sum_{\mathbf{i=1}}^{\mathbf{i=n-1}}x_{\mathbf{i}} + (6+5n-n^2)x_{\mathbf{n}}.
\end{equation}

\noindent This ellipsoid is unique.
\end{theorem}

Michael Greene, a Harvard undergraduate who worked in the summer of 2001 with Rybnikov
in the \textbf{R}esearch
\textbf{E}xperience for \textbf{U}ndergraduates Project at Cornell  proved that
 $\mathbf{x^{t} Q_7 x+Lx=0}$ is Delaunay by checking it directly with a computer program that he wrote in
\textsc{Mathematica} system (this program is available by request from Rybnikov:
krybniko@cs.uml.edu). Recently, we proved the Delaunay property for all of the series.

In  2002 Bob Erdahl has found an important realization of $\Upsilon _{n}$ as a section
of a centrally-symmetric polytope $C_n$ in $\mathbb{Z}^{n+1}.$ For $n=6,$ this
realization is, in fact, $2_{21}$ (Gosset 6D polytope) as a section of $ 3_{21}$
(Gosset 7D polytope).

In general, consider $\mathbb{Z}^{n}$ as a subspace of $\mathbb{Z}^{n+1}$ defined by $%
x_{n+1}=1.$ Then take two copies of the original $\Upsilon_{n}:$ a shift of $\Upsilon
_{n}$ by vector $[0^{n};1]$ and an inverted copy of $ \Upsilon_{n}$, i.e.
$-\Upsilon_{n}.$ It turns out that one can complement these two polytopes with two
points, so that the convex hull of the resulting set is a perfect $(n+1)$-dimensional
Delaunay polytope.

The two points that need to be added to the vertices of $\Upsilon _{n}+[0^{n};1]$ and
$-\Upsilon _{n}$ are $[1^n;(n-3)]$ and $ [-1^n;-(n-4)].$ Denote the resulting point set
by $\vert C_{n+1}.$ Points of $\vert C_{n+1}$ \ can be partitioned into pairs
corresponding to the segments passing through $[0^{n};\frac{1}{2}]$. These segments are
the diagonals of the $(n+1)$-dimensional polytope $C_{n+1}$, mentioned above. Let us
now consider  lattice $\Lambda_{n+1}$ generated by the diagonals of $C_{n+1}. $
$2\mathbb{ Z}^{n+1}$ is a sublattice of $\Lambda_{n+1} :$ if one can prove that the
diagonals are the minimal vectors of $\Lambda_{n+1} $ with respect to some quadratic
form $\varphi$, then one has a proof that they are minimal vectors in their parity
class in $\mathbb{Z}^{n+1}$ relative to this form $\varphi$. Vectors $ B=\{[0^{n};1],$
$[2,0^{n-1};1]$, $[2^{n};2n-7]\}$ are linearly independent and, therefore, form a
rational basis of $\Lambda_{n+1} .$ The Gramm
matrix of $B,$ with respect to $Q_n(\x',\x') \bigoplus x_{n+1}^2$ where $\x=[\x';x_{n+1}]$, has an astonishingly simple form: it is $%
I+\alpha J,$ where $\alpha =\frac{1+\binom{n-4}{2}}{8(n-5)}$ and $J$ is the $n+1$ by
$n+1$ matrix of all 1's.

\bigskip

The philosophy behind the construction of supertopes is in generalizing the affine
structure of $2_{21}$ to a general dimension. While most important affine properties
are preserved in this generalization, the group-theoretic ones are not.\ In our
research on supertopes  we have been using software \textsc{Polymake} written by
Gawrilow and Joswig in TU Berlin.

We hope that one can simplify the analysis of supertopes by using the symmetries of the
supertopes, in particular, the so-called ''long
triangles'' (e.g $[1^{2};0^{n-3};-1]$, $[0^{2};1^{2};0^{n-5};-1]$, $%
[0^{4};1^2;0^{n-7};-1]$ in $\Upsilon_n$). These triangles are similar to 45 long
equatorial
triangles (of norm 2) on the boundary of $2_{21}$, the Delaunay cell of $%
E_{6}$ and an affine image of our 6-supertope. Right now we are checking if these
triangles are playing in an $n$-supertope a role, similar to that of the long triangles
of $2_{21}$: the projection of $2_{21}$ on the 4-plane perpendicular to the long
triangle is a 24-cell (Ivi\v{c} Weiss, 1986). In the summer of 2002, within the
framework of Cornell REU,\ Rybnikov worked with Cornell undergraduate student Joseph
Palin on connections between "long triangles" and  arithmetic properties of
inhomogeneous quadratic forms defined by polytopes $\Upsilon _{n}$.

Long triangles (of norm 2) in $E_{n}$ are interesting for another reason.
Arrangements of equilateral triangles of side length $\sqrt{2}$ are crucial
in Blichfeldt's (1935) argument (published without proof) that the densest
lattices for $n=6,7,8$ must have $D_{5}$ as a sublattice of the same minimum
(see Vetchinkin (1982) for a proof of Blichfeldt's result).

Voronoi vectors of a positive definite form $\varphi (\mathbf{x})$ are, by definition,
the minimal integral vectors, relative to $\varphi (\mathbf{x})$, in their parity
classes, i.e. $\textrm{mod}\;2$. Roughly speaking, each Delaunay face with inner
diagonals (such faces are called \emph{primitive elements} of the tiling) is the convex
hull of Voronoi vectors from the same parity class with common midpoint; these
diagonals are the only minimal vectors in their parity class. The converse is also
true. Although these facts have been first observed and used by Voronoi (1909),
Baranovskii (1991) was the first to formulate and prove them as a theorem. The (partial
for $n=5$)
classification of $L$-types for $n<6$ was obtained via the theory of \emph{%
primitive elements}.  We are trying to generalize the Voronoi-Baranovskii theory of
parity classes to \emph{triangles} of Voronoi vectors. Voronoi triangles seem to be
very important for the problem of densest lattice packings (methods of Blichfeldt and
Vetchinkin) and the theory of perfect forms in general. The problem here is that not
all Voronoi triangles from can play a role similar to that of Voronoi vectors. $ E_{6}$
and the series of forms from the above theorem, basically, serves as a laboratory for
this project.

We hope that our work on positive inhomogeneous quadratic forms may lead to a better
understanding of the structure of lattices and the relationship between the theories of
$L$-types and perfect forms. This could shed new light on how the packing and covering
problems are related to each other.


\begin{thebibliography}{99}

\bibitem{Alex2001} V.~Alexeev (2001), On Extra Components in the Toroidal Compactification of
$A_g$. \emph{Moduli of Abelian Varieties (Texel island)}, Progress in Mathematics
\textbf{195,} Birkhauser, Boston, pp. 1--7. Also on arXiv.org: math.AG/9905142.

\bibitem{Alex2002} V.~Alexeev (2002), Complete moduli in the presence of semiabelian group
action. \emph{Annals of Mathematics}  \textbf{155}, no. 3, pp. 611--708.  Also on
arXiv.org: math.AG/9905103.

\bibitem{Anz}  M.M. Anzin (1992), On Variations of Positive Quadratic Forms (with
applications to the study of perfect forms), \textit{Proc. Steklov Inst. Math. }
\textbf{196,} no. 4, pp. 11--27.

\bibitem{Baran91-92}  E.~P.~Baranovskii (1991), Partition of Euclidean spaces into $L$
-polytopes of certain perfect lattices. (Russian) Discrete geometry and topology
(Russian). \emph{Trudy Mat. Inst. Steklov. } \textbf{196}, pp. 27--46. Translated in
\emph{Proc. Steklov Inst. Math. } \textbf{196,} (1992), no. 4, pp. 29--51.

\bibitem{Barn1957}  E.~S.~Barnes (1957), The complete enumeration of extreme senary forms.
\textit{Philos. Trans. Roy. Soc. London Ser. A} \textbf{249}, pp. 461--506.

\bibitem{BarnDick1967}  E.~S.~Barnes, T.~J.~Dickson (1967), Extreme coverings of $n$-space by
spheres. \textit{J. Austral. Math. Soc.} \textbf{7}, pp. 115--127.

\bibitem{Blich1935}  H.~F.~Blichfeldt (1935), The minumum values of positive quadratic forms
in six, seven and eight variables. \textit{Math. Z}. \textbf{39}, pp. 1--15.

\bibitem{CSperfect1988}  J.~H.~Conway, N.~J.~A.~Sloane (1988), Low-dimensional lattices. III.
Perfect forms.
 \textit{Proc. Roy. Soc. London Ser. A} \textbf{418}, no. 1854, pp. 43--80.

\bibitem{CSthird1999}  J.~H.~Conway, N.~J.~A.~Sloane (1999), \textit{\ Sphere packings,
lattices and groups.} Third edition. Grundlehren der Mathematischen Wissenschaften
\textbf{290,} Springer-Verlag, New York.

\bibitem{CoxAnnRed1934-95}  H.~S.~M.~Coxeter (1934), Discrete Groups Generated by Reflections.
\emph{Annals of Math.} \textbf{35}, pp. 588--621. Repr.\ in \emph{Kaleidoscopes:
Selected Writings of H.~S.~M.~Coxeter}, Eds.: F.~A.~Sherk \emph{et al.}, Wiley, New
York, 1995.

\bibitem{Cox1951EF}  H.~S.~M.~Coxeter (1951), Extreme forms. \textit{\ Canadian J. Math.}
no. 3, pp. 391--441.

\bibitem{Cox1973}  H.~S.~M.~Coxeter (1973), \textit{Regular Polytopes.} 3rd
edition. Dover, New York.

\bibitem{DELTor1924-28}  B.~N.~Delaunay [Delone] (1924), Sur la sph\`{e}re vide, in: \textit{
Proceedings of the International Congress of Mathematicians, Toronto, 1924,} University
of Toronto Press, Toronto, (1928), pp. 695--700.

\bibitem{DEboth1937}  B.~N.~Delaunay [Delone] (1937), The geometry of positive quadratic
forms, \textit{Uspekhi Mat. Nauk} \textbf{3} (1937), pp. 16--62 and \textbf{4} (1938),
pp. 102-164.

\bibitem{DRfour1963R}  B.~N.~Delone, S.~S.~Ry\v{s}kov (1963), Solution of the problem on the
least dense lattice covering of a 4-dimensional space by equal spheres. (Russian)
\emph{Dokl. Akad. Nauk SSSR} \textbf{152,} pp. 523--524.

\bibitem{DDRS1970E}  B.~N.~Delone, N.~P.~Dolbilin,  S.~S.~Ry\v{s}kov, M.~I.~\v{S}togrin (1970),
 A new construction of the theory of lattice coverings of an $n$-dimensional space by
 congruent balls. (Russian) Izv. Akad. Nauk SSSR Ser. Mat. \textbf{34.}

\bibitem{DRtwo1971} B.~N.~Delone, S.~S.~Ry\v{s}kov (1971),  Extremal problems of the theory of
positive quadratic forms. \textit{Trudy Mat. Inst. Steklov} \textbf{112}, pp. 203--223,
387.

\bibitem{DG1993}  M.~Deza, V.~P.~Grishukhin (1993), Hypermetric graphs. \emph{Quart.
J. Math. Oxford Ser. (2)} \textbf{44}, no. 176, pp. 399--433.

\bibitem{DG1995}  M.~Deza, V.~P.~Grishukhin (1995), Delaunay polytopes of cut
lattices. \emph{Linear Algebra Appl.} \textbf{226/228,} pp. 667--685.

\bibitem{DG1994-96}  M.~Deza, V.~P.~Grishukhin (1996), Cut lattices and equiangular
lines. Discrete metric spaces (Bielefeld, 1994). \textit{European J. Combin.}
\textbf{17} (1996), no. 2-3, pp. 143--156.


\bibitem{DGL1991-92}  M.~Deza, V.~P.~Grishukhin, M.~Laurent (1992), Extreme hypermetrics
and $L$-polytopes. Sets, graphs and numbers (Budapest, 1991), 157--209, \textit{Colloq.
Math. Soc. Janos Bolyai}, \textbf{60}, North-Holland, Amsterdam, 1992.

\bibitem{DGL1992-3-5}  M.~Deza, V.~P.~Grishukhin, M.~Laurent (1995), Hypermetrics in
geometry of numbers. \emph{\ Combinatorial optimization (New Brunswick, NJ,
1992--1993)}, pp. 1--109, \emph{DIMACS Ser. Discrete Math. Theoret. Comput. Sci.}
\textbf{20}, Amer. Math. Soc., Providence, RI, 1995.

\bibitem{DL1997}  M.~Deza, M.~Laurent (1997), \textit{Geometry of cuts and metrics.}
Algorithms and Combinatorics \textbf{15.} Springer-Verlag, Berlin.

\bibitem{Dick1968}  T.~J.~Dickson (1968), A sufficient condition for an extreme covering
of $n$-space by spheres. \textit{J. Austral. Math. Soc.} no. 8, pp. 56--62.

\bibitem{Dick1972}  T.~J.~Dickson (1972), On Voronoi reduction of positive definite
quadratic forms. \textit{\ J. Number Theory} no. 4, pp. 330--341.

\bibitem{Diet1975}  U.~Dieter (1975), How to calculate shortest vectors in a lattice.
\textit{Mathematics of Computation}, \textbf{29}, no. 131, pp. 827--833.

\bibitem{E1975}  R. Erdahl (1975), \emph{A convex set of second-order inhomogeneous
polynomials with applications to quantum mechanical many body theory.} Mathematical
Preprint \#1975-40, Queen's University, Kingston, Ontario.

\bibitem{E1992}  R.~Erdahl (1992), A cone of inhomogeneous second-order polynomials,
\textit{Discrete Comput. Geom.} \textbf{8}, no. 4, pp. 387--416.

\bibitem{E2000}  R.~M.~Erdahl (2000), \textit{A structure theorem for Voronoi
polytopes of lattices}, Talk at the sectional meeting of the AMS, Session on Discrete
Geometry, Toronto, September 22--24.

\bibitem{E2001}  R.~M.~Erdahl (2001), \textit{On the tame facet of the perfect domain
}$E_{6}^{\ast }$, Plenary talk at the Seventieth Birthday Conference for Sergei
Ryshkov, Steklov Institute of Mathematics, Moscow, January 24--27.

\bibitem{ER-RendicPal2002} R.~M.~Erdahl, K.~Rybnikov (2002), Voronoi's Hypothesis on Perfect Forms and L-types.
\emph{Rendiconti del Circolo  Matematiko di Palermo}, Serie II, Tomo LII, part I, pp.
279--296.  Full-length version available on arXiv.org: \textit{math.MG/0112098}.

\bibitem{ERProc2002E} R.~M.~Erdahl, K.~Rybnikov (2002b),  New Infinite Series of Perfect Quadratic Forms
and Big Delaunay Simplexes in $\mathbb{Z}^{n}$. \textit{Proc. of Steklov Institute of
Math.}  \textbf{239}, no. 4, pp. 159--167. Also on arXiv.org: \textit{math.MG/0112098}.

\bibitem{HaaZ2000}  C.~Haase, G.~Ziegler (2000), On the maximal width of empty lattice
simplices. Combinatorics of polytopes. \textit{European J. Combin.} \textbf{21}, no. 1,
pp. 111--119.

\bibitem{KZ1873}  A.~Korkine, G.~Zolotareff (1873), Sur les formes quadratiques,
\textit{Math. Ann.} 6, pp. 366--389.

\bibitem{LemmSeid1973}  P.~H.~W.~Lemmens, J.~J.~Seidel (1973), Equiangular lines. \emph{J.
Algebra} 24, pp. 494--512.

\bibitem{Mart1996F}  J.~Martinet (1996), \textit{Les Reseaux Parfaits des espaces Euclidiens.}
Masson, Paris.

\bibitem{Nie1973}  H-V.~Niemeier (1973), Definite quadratische Formen der Dimension 24 und
Diskriminate 1. (German. English summary) \textit{J. Number Theory} no. 5, pp.
142--178.

\bibitem{Rgroups1972} S.~S.~Ryshkov (1972), Maximal finite groups of integral n x n
matrices  and full groups of integral automorphisms of positive quadratic forms
(Bravais models). \emph{Trudy Mat. Inst. Steklov} Vol. \textbf{128,} pp. 183--211 (in
Russian). English translation in \emph{Proc. Steklov Inst. Math.} (1972) Vol.
\textbf{128,} pp. 217--250.

\bibitem{RMink1973-76}  S.~S.~Ryshkov (1973), The perfect form $A_{n}^{k}:$ existence of
lattices with a nonfundamental partition simplex; existence of perfect forms that are
not reducible in the sense of Minkowski to a form with identical diagonal coefficients.
(English translation of Russian, 1973, original) \emph{J. Soviet Math.} 6 (1976), pp.
672--676.

\bibitem{R1999}  S.~S.~Ryshkov (1999), Direct geometric description of $n$%
-dimensional Voronoi parallelohedra of the second type, (Russian) \textit{ Uspekhi Mat.
Nauk} \textbf{54}, no. 1 (325), pp. 263--264; English translation in \textit{Russian
Math. Surveys} \textbf{54} (1999), no. 1, pp. 264--265.

\bibitem{R1998}  S.~S.~Ryshkov (1998), On the structure of a primitive
parallelohedron and Voronoi's last problem, (Russian) \textit{Uspekhi Mat. Nauk}
\textbf{53}, no. 2 (320), pp. 161--162; translation in \textit{ Russian Math. Surveys}
\textbf{53} (1998), no. 2, pp. 403--405.

\bibitem{Sen1995}  M.~Senechal (1995), \textit{Quasicrystals and geometry. }Cambridge
University Press, Cambridge, UK.

\bibitem{Vetch}  N.~M.~Vet\v{c}inkin  (1982), Uniqueness of classes of positive
quadratic forms, on which values of Hermite constants are reached for $\ 6\leq n\leq
8.$ \textit{Proceeding of the Steklov Institute of Math.} \textbf{152}, No. 3, pp.
37--95.

\bibitem{VorDeux1908-09}  G.~F.~Voronoi (1908), Nouvelles applications des param\`{e}ters
continus \`{a} la th\'{e}orie des formes quadratiques, Deuxi\`{e}me memoire, \textit{J.
Reine Angew. Math.} \textbf{134} (1908), pp. 198--287 and \textbf{136} (1909), pp.
67--178.

\bibitem{VorThree1952}  G.~F.~Voronoi (1952), \textit{Sobranie socinenii v treh tomah}
[Collected works in three volumes] \textbf{Vol. 2.} Introduction and notes by
B.~N.~Delaunay. (all in Russian), Kiev.

\bibitem{AI-W1986}  A.~Ivi\v{c}~Weiss (1986), A four-dimensional projection of the polytope
$2_{21}.$ \textit{C. R. Math. Rep. Acad. Sci. Canada} \textbf{8}, no. 6, pp. 405--410.
\end{thebibliography}
\end{document}